\documentclass[reqno,report,12pt]{amsart} 
\usepackage{mathrsfs}
\usepackage[margin=1in]{geometry}
\usepackage{amssymb,amsmath,amsthm,amsfonts}
\usepackage{latexsym}
\usepackage{amssymb}
\usepackage{amsmath}
\usepackage{amsfonts}
\usepackage{amsthm}
\usepackage{amscd}
\usepackage{texdraw}
\usepackage{color}
\usepackage{cite}
\usepackage[centertags]{amsmath}
\usepackage[colorlinks=true,urlcolor=blue,linkcolor=red,citecolor=red]{hyperref}
\DeclareMathOperator{\ima}{Im}
\newlength{\defbaselineskip}
\newcommand{\setlinespacing}[1]%
           {\setlength{\baselineskip}{#1 \defbaselineskip}}

\theoremstyle{plain}
\newtheorem{theorem}{Theorem}
\newtheorem{lemma}{Lemma}

\newtheorem{corollary}{Corollary}

\newtheorem{remark}{Remark}

{\allowdisplaybreaks
\begin{document}
\title[Weighted inversion of vector valued Dirichlet series]{Weighted inversion of vector valued Dirichlet series}

\author[P. A. Dabhi]{Prakash A. Dabhi}
\address{Department of Basic Sciences, Institute of Infrastructure, Technology, Research And Management (IITRAM), Maninagar (East), Ahmedabad - 380026, Gujarat, India}
\email{lightatinfinite@gmail.com, prakashdabhi@iitram.ac.in}

\author[K. B. Solanki]{Karishman B. Solanki}
\address{Department of Mathematics, Indian Institute of Technology Ropar, Rupnagar, Punjab - 140001, India}
\email{karishsolanki002@gmail.com, staff.karishman.solanki@iitrpr.ac.in}

\thanks{\textit{Corresponding author}: Karishman B. Solanki}



\date{}

\dedicatory{}

\commby{}

\begin{abstract}
Let $\Lambda\subset[0,\infty)$ be an additive semigroup with $0\in\Lambda$, $\omega$ be an admissible weight on $\Lambda$, $\mathcal A$ be a unital Banach algebra, and let $f(s)=\sum_{\lambda\in\Lambda} f_\lambda e^{-\lambda s}$ for $s\in\mathcal{H}=\{j+it\in\mathbb{C}:j\geq0\}$ be a generalized Dirichlet series satisfying $\|f\|_\omega=\sum_{\lambda\in\Lambda}\|f_\lambda\|\omega(\lambda)<\infty,$ where $f_\lambda\in\mathcal{A}$ for all $\lambda\in\Lambda$. We take $\mathcal{A}$ to be a commutative complex Banach algebra (with $\Lambda=\log\mathbb{N}$) and $M_d(\mathcal{X})$ - the Banach algebra of $d \times d$ matrices having entries from $\mathcal{X}$, where $\mathcal{X}$ is either the complex plane or the real algebra of bicomplex numbers or quaternions, and show that $f$ is invertible if and only if the closure of the image of $f$ is contained in the set of all invertible elements of $\mathcal{A}$.
\end{abstract}

\maketitle


\section{Introduction} \label{sec:intro}

Let $\Lambda\subset[0,\infty)$ be an additive semigroup with $0\in\Lambda$, $\mathcal{H}=\{s=j+it:j\geq0\}$ be the right half plane of the complex plane $\mathbb C$, and let $(\mathcal{A},\|\cdot\|)$ be a unital (complex or real) Banach algebra with unit $1_\mathcal{A}$. Let $\mathcal{D}(\Lambda,\mathcal{A})$ be the collection of all functions $f:\mathcal{H}\to\mathcal{A}$ having the generalized Dirichlet series form 
\begin{align} \label{eqn:dirform}
    f(s)=\sum_{\lambda\in\Lambda} f_\lambda e^{-\lambda s} \quad (s\in\mathcal{H}) \quad \text{satisfying} \quad \|f\|_1=\sum_{\lambda\in\Lambda}\|f_\lambda\|<\infty,
\end{align}
where $f_\lambda\in\mathcal{A}$ for all $\lambda\in\Lambda$. Then $\mathcal{D}(\Lambda,\mathcal{A})$ is a unital Banach algebra with the norm $\|\cdot\|_1$ and the pointwise multiplication. Observe that if $\displaystyle f(s)=\sum_{\lambda \in \Lambda}f_\lambda e^{-\lambda s}$ and $g(s)=\displaystyle \sum_{\lambda \in \Lambda}g_\lambda e^{-\lambda s}$ are in $\mathcal D(\Lambda,\mathcal A)$, then 
\begin{align} \label{eq:pintwiseproduct}
    (fg)(s)=f(s)g(s)=\sum_{\lambda \in \Lambda}\left(\sum_{\substack{\lambda_1,\lambda_2 \in \Lambda \\ \lambda_1+\lambda_2=\lambda}} f_{\lambda_1} g_{\lambda_2}\right)e^{-\lambda s} \quad (s\in \mathcal H).
\end{align}

Let $\ell^1(\Lambda,\mathcal A)=\{ (f_\lambda)_{\lambda\in\Lambda} : f_\lambda \in \mathcal A \text{ for all } \lambda\in\Lambda \, \text{and} \, \|f\|_1=\sum_{\lambda\in \Lambda}\|f_\lambda\|<\infty\}$. Then $\ell^1(\Lambda,\mathcal A)$ is a unital Banach algebra with the above norm and the convolution product given by
\begin{align} \label{eq:convolutionproduct}
(f\star g)_\lambda = \sum_{\substack{\lambda_1,\lambda_2 \in \Lambda \\ \lambda_1+\lambda_2=\lambda}} f_{\lambda_1} g_{\lambda_2}\quad(f,g \in \ell^1(\Lambda,\mathcal A)).    
\end{align}
The Banach algebras $\ell^1(\Lambda,\mathcal A)$ and $\mathcal D(\Lambda,\mathcal A)$ are isometrically isomorphic via the map 
\begin{align*}
\ell^1(\Lambda,\mathcal{A})\ni (f_\lambda)\mapsto f(s)=\sum_{\lambda \in \Lambda}f_\lambda e^{-\lambda s}\in \mathcal D(\Lambda,\mathcal A).
\end{align*}

\begin{remark}
    The fact $0\in\Lambda$ ensures 
    \begin{enumerate}
        \item[(a)] the products in \eqref{eq:pintwiseproduct} and \eqref{eq:convolutionproduct} are well defined, and 
        \item[(b)] the existence of the unit element in $\mathcal{D}(\Lambda,\mathcal{A})$ and $\ell^1(\Lambda,\mathcal A)$ which are given by $f(s)=1_\mathcal{A}$ $(s\in\mathcal{H})$ and $\delta_0(0)=1_\mathcal{A}$ with zero otherwise, respectively.
    \end{enumerate}
\end{remark}

Consider $\mathbb{N}$ as a semigroup with multiplication and the Banach algebra $\ell^1(\mathbb{N},\mathcal{A})$ with the multiplication 
\begin{align*}
    (a\star b)(n)=\sum_{\substack{k,l\in\mathbb{N} \\ kl=n}} a_kb_l \quad (a=(a_n)_{n\in\mathbb{N}}, b=(b_n)_{n\in\mathbb{N}} \in \ell^1(\mathbb{N},\mathcal{A})).
\end{align*} 
Take $\Lambda=\log\mathbb{N}$. Then $\mathcal{D}(\log\mathbb{N},\mathcal{A})$ is the Banach algebra of ordinary Dirichlet series, that is, $f\in\mathcal{D}(\log\mathbb{N},\mathcal{A})$ is of the form 
$$f(s)=\sum_{n\in\mathbb{N}} f_n n^{-s} \quad (s\in\mathcal{H}),$$ 
and the multiplication of $f,g\in\mathcal{D}(\log\mathbb{N},\mathcal{A})$ becomes 
\begin{align*}
    f(s)g(s)=\sum_{n\in\mathbb{N}} \left( \sum_{\substack{k,l\in\mathbb{N}\\ kl=n}} f_k g_l \right) n^{-s} \quad (s\in\mathcal{H})
\end{align*}
For notational convenience and keeping up with literature, we shall write $\mathcal{D}(\mathbb{N},\mathcal{A})$ instead of $\mathcal{D}(\log\mathbb{N},\mathcal{A})$ as it becomes isometrically isomorphic to $\ell^1(\mathbb{N},\mathcal{A})$ defined above.

In \cite{he}, Hewitt and Williamson proved that $f\in\mathcal{D}(\mathbb{N},\mathbb{C})$ is invertible if and only if $f$ is bounded away from zero in the absolute value. An elementary proof of it is obtained in \cite{gn}. The following is our first main theorem, which is a generalization of the above theorem of Hewitt and Williamson \cite[Theorem 1]{he} for a commutative unital complex Banach algebra $\mathcal{A}$. By $G(\mathcal{A})$ denote the open subset of $\mathcal{A}$ containing all invertible elements of $\mathcal{A}$.

\begin{theorem}\label{commutative}
Let $\mathcal{A}$ be a commutative unital complex Banach algebra, $a=(a_n)_{n\in\mathbb{N}} \in \ell^1(\mathbb N,\mathcal A)$, and let $f_a(s)=\sum_{n\in \mathbb N} a_n n^{-s}$ for all $s \in \mathcal H$. Then the following statements are equivalent.
\begin{enumerate}
\item $a$ is invertible in $\ell^1(\mathbb{N},\mathcal{A})$.
\item $\overline{\ima(f_a)}=\overline{\{f_a(s):s\in\mathcal{H}\}}\subset G(\mathcal{A})$, where the overline denotes the closure in $\mathcal A$.
\item There is $(b_n)_{n\in\mathbb{N}} \in \ell^1(\mathbb N,\mathcal A)$ such that $(f_a(s))^{-1}=\sum_{n\in\mathbb{N}} b_n n^{-s}$ for all $s\in\mathcal{H}$.
\end{enumerate}
In particular, a function $f\in\mathcal{D}(\mathbb{N},\mathcal{A})$ is invertible if and only if $\overline{\ima(f)}\subset G(\mathcal{A})$.    
\end{theorem}

In \cite{edw}, Edwards generalized the result of Hewitt and Williamson for arbitrary $\Lambda$ again for $\mathcal{A}=\mathbb{C}$ as stated below.

\begin{theorem}[Edward's theorem]\cite{edw}\label{th1}
An element $f\in\mathcal{D}(\Lambda,\mathbb{C})$ has an inverse in $\mathcal{D}(\Lambda,\mathbb{C})$ if and only if $f$ is bounded away from $0$ in absolute value.
\end{theorem}

We do not know whether Theorem \ref{commutative} is true for arbitrary $\Lambda$ and a commutative Banach algebra $\mathcal{A}$.

Now, we turn towards the weighted versions. A \emph{\textbf{weight}} on $\Lambda$ is a map $\omega:\Lambda\to[1,\infty)$ satisfying 
\begin{align*}
\omega(0)=1 \quad \text{and} \quad \omega(\lambda_1+\lambda_2) \leq \omega(\lambda_1) \omega(\lambda_2) \quad \text{for all} \quad \lambda_1,\lambda_2\in\Lambda.
\end{align*} 
A weight $\omega$ on $\Lambda$ is an \emph{\textbf{admissible weight}} if $\displaystyle \lim_{n\to\infty} \omega(n\lambda)^\frac{1}{n}=1$ for all $\lambda\in\Lambda$. For an admissible weight $\omega$ on $\Lambda$, let $\mathcal{D}(\Lambda,\omega,\mathcal{A})$ be the collection of all function $f:\mathcal{H}\to\mathcal{A}$ of the form $$f(s)=\sum_{\lambda\in\Lambda} f_\lambda e^{-\lambda s} \quad (s\in\mathcal{H}) \quad \text{satisfying} \quad \|f\|_\omega=\sum_{\lambda\in\Lambda}\|f_\lambda\| \omega(\lambda)<\infty.$$
Then $\mathcal D(\Lambda,\omega,\mathcal A)$ is a subalgebra of $\mathcal D(\Lambda,\mathcal A)$ and $\mathcal D(\Lambda,\omega,\mathcal A)$ is a Banach algebra with the (weighted) norm $\|\cdot\|_\omega$.
Gl\"ockner and Lucht \cite{glo} proved the following weighted version of Edward's theorem for an admissible weight. Edward's theorem follows by taking $\omega$ to be the trivial weight $\omega(\lambda)=1$ for all $\lambda\in\Lambda$.

\begin{theorem}\cite{glo} \label{th2}
Let $\Lambda\subset[0,\infty)$ be an additive semigroup with $0\in\Lambda$, and let $\omega$ be an admissible weight on $\Lambda$. A function $f\in\mathcal{D}(\Lambda,\omega,\mathbb{C})$ has an inverse in $\mathcal{D}(\Lambda,\omega,\mathbb{C})$ if and only if $f$ is bounded away from $0$ in absolute value.
\end{theorem}

Let $\mathcal{A}$ be a unital subalgebra of a unital Banach algebra $\mathcal{B}$ both having the same unit. We say $\mathcal{A}$ is \textit{\textbf{inverse-closed}} or \textbf{\textit{spectrally invariant}} in $\mathcal{B}$ if $a\in\mathcal{A}$ is invertible in $\mathcal{B}$ implies that $a$ is invertible in $\mathcal{A}$. So, combining Theorem \ref{th1} and Theorem \ref{th2}, it follows that $\mathcal{D}(\Lambda,\omega,\mathbb{C})$ is inverse-closed in $\mathcal{D}(\Lambda,\mathbb{C})$ provided $\omega$ is an admissible weight on $\Lambda$. The inverse-closedness of Banach algebras has been studied extensively and are of importance in the study spectral theory linear and non-linear operators on $C^\star$-algebras and Banach algebras in general. We refer the reader to \cite{barnes, beaver, folland, Gro2, KB3, shin}, and reference therein.

We shall prove some vector valued analogues of Theorem \ref{th2} for some classes of unital Banach algebras. First one is by taking $\mathcal{A}$ to be a commutative unital complex Banach algebra with $\Lambda=\log\mathbb{N}$ and it is weighted version of Theorem \ref{commutative}. Note that the non-weighted results will simply follow by taking the trivial weight $\omega\equiv1$ as it is an admissible weight. 

\begin{theorem}\label{weighted commutative}
Let $\omega$ be an admissible weight on $\mathbb{N}$, let $\mathcal{A}$ be a commutative unital complex Banach algebra, and let  $f\in\mathcal{D}(\mathbb{N},\omega,\mathcal{A})$. Then $f$ is invertible in $\mathcal{D}(\mathbb{N},\omega,\mathcal{A})$ if and only if $\overline{\ima(f)}\subset G(\mathcal{A})$.
\end{theorem}

We have an immediate corollary as stated below.
\begin{corollary}
Let $\omega$ be an admissible weight on $\mathbb{N}$, and let $\mathcal{A}$ be a commutative unital complex Banach algebra. Then $\mathcal{D}(\mathbb{N},\omega,\mathcal{A})$ is inverse-closed in $\mathcal{D}(\mathbb{N},\mathcal{A})$.
\end{corollary}

As stated earlier, we do not know whether the above Theorem \ref{weighted commutative} and hence Theorem \ref{commutative} is true for arbitrary $\Lambda$ and noncommutative $\mathcal{A}$. But we shall deal with three noncommutative Banach algebras here and get an affirmative answer. Out of these three, one is a complex algebra while two are real algebras. And this open the questions whether these results follows when $\mathcal{A}$ is a real algebra. We do not know the answer to it also as the real algebras here are dealt with the ``\textit{slice}" approach which eventually lands us back in the complex algebra setting.

Now, moving to the noncommutative algebra analogues, the first case is for complex matrices. For $d\in\mathbb{N}$, take $\mathcal{A}=M_d(\mathbb{C})$, the set of all $d\times d$ complex matrices. Then $\mathcal{A}$ is a noncommutative unital complex Banach algebra with operator norm and the $d\times d$ identity matrix $I_d$ as the unit. We have the following result for it.

\begin{theorem}\label{matrix}
Let $d\in\mathbb{N}$, $\Lambda\subset[0,\infty)$ be an additive semigroup with $0\in\Lambda$, and let $\omega$ be an admissible weight on $\Lambda$. An element $f\in\mathcal{D}(\Lambda,\omega,M_d(\mathbb{C}))$ has an inverse in $\mathcal{D}(\Lambda,\omega,M_d(\mathbb{C}))$ if and only if $\overline{\ima(f)}=\overline{\{f(s):s\in\mathcal{H}\}}\subset GL_d(\mathbb{C})$, the collection of all invertible matrices from $M_d(\mathbb{C})$.
\end{theorem}

Next, we deal with the four dimensional real algebra of bicomplex numbers $\mathbb{BC}$. It is generated by two commutating imaginary units $i$ and $j$ along with the \emph{\textbf{hyperbolic unit}} $k:=ij=ji$. That is an element $Z\in\mathbb{BC}$ is written as 
\begin{align*}
    Z=x_1 + x_2 i +x_3 j + x_4 k, \quad \text{where} \quad x_1,x_2,x_3,x_4\in\mathbb{R}.
\end{align*}
The element $k$ is called a hyperbolic unit because $k^2=1$.

Let $\mathbf{e}_1=\frac{1+k}{2}$ and $\mathbf{e}_2=\frac{1-k}{2}$ be the elements of $\mathbb{BC}$. Then $\mathbf{e}_1\cdot\mathbf{e}_2=0$, $\mathbf{e}_1^2=\mathbf{e}_1$, $\mathbf{e}_2^2=\mathbf{e}_2$, $\mathbf{e}_1+\mathbf{e}_2=1$ and $\mathbf{e}_1-\mathbf{e}_2=k$. The first equality implies that $\mathbb{BC}$ is not a division algebra. 
The set $\{\mathbf{e}_1,\mathbf{e}_2\}$ forms a basis of $\mathbb{BC}$ over the complex plane $\mathbb{C}_i=\{x+iy:x,y\in\mathbb{R}\}$ and it is known as \emph{\textbf{idempotent basis}}. The set $\{1,j\}$ also forms a basis of bicomplex numbers over $\mathbb{C}_i$. The dimension of $\mathbb{BC}$ over $\mathbb{C}_i$ is two. We shall just write $\mathbb{C}$ instead of $\mathbb{C}_i$ whenever $i$ in consideration is clear.
The \emph{\textbf{idempotent representation}} of a bicomplex number $Z=z_1+ j z_2$ $(z_1,z_2\in\mathbb{C})$ is given by 
\begin{align*}
    Z=\lambda_1 \mathbf{e}_1 + \lambda_2 \mathbf{e}_2, \quad \text{where} \quad \lambda_1=z_1- i z_2 \quad \text{and} \quad \lambda_2=z_1 + i z_2.
\end{align*} It follows that $z_1=\frac{\lambda_1+\lambda_2}{2}$  and $z_2=\frac{i(\lambda_1-\lambda_2)}{2}$.

The reason for choosing this particular representation of bicomplex numbers is that the addition and multiplication of bicomplex numbers can be realized componentwise in the idempotent representation. It means that for $Z=\lambda_1 \mathbf{e}_1 + \lambda_2 \mathbf{e}_2$ and $W=\mu_1 \mathbf{e}_1 + \mu_2 \mathbf{e}_2$ in $\mathbb{BC}$ with $\lambda_1,\lambda_2,\mu_1,\mu_2\in\mathbb{C}$ and $n \in \mathbb N$, we have 
\begin{align} \label{eq:BCrelations} \nonumber
    &Z+W=(\lambda_1+\mu_1)\mathbf{e}_1 + (\lambda_2+\mu_2)\mathbf{e}_2, \\ \nonumber 
    &Z\cdot W=(\lambda_1 \mu_1)\mathbf{e}_1 + (\lambda_2 \mu_2)\mathbf{e}_2 \quad \text{and} \\ & Z^n=\lambda_1^n\mathbf{e}_1 + \lambda_2^n \mathbf{e}_2.
\end{align}
Also, if $Z$ is an invertible bicomplex number, that is, $\lambda_1 \lambda_2\neq 0$, then $Z^{-1}=\lambda_1^{-1}\mathbf{e}_1 + \lambda_2^{-1}\mathbf{e}_2$. The norm on $\mathbb{BC}$ that we shall consider here is the \textit{\textbf{dual Lie norm}} given by $\|Z\|=|\lambda_1|+|\lambda_2|$.

Let $d\in\mathbb{N}$, $\omega$ be a weight on $\Lambda$, and let $\mathcal{A}=M_d(\mathbb{BC})$ be the collection of $d\times d$ matrices having entries from $\mathbb{BC}$. Given a matrix $M\in M_d(\mathbb{BC})$, we can write its idempotent representation as 
\begin{align*}
    M=M_1\mathbf{e}_1 + M_2\mathbf{e}_2, \text{ where } M_1,M_2 \in M_d(\mathbb{C}_i).
\end{align*}
It can be easily verified that the relations in \eqref{eq:BCrelations} also follows for the elements of $M_d(\mathbb{BC})$ and the norm of $M$ is given by $\|M\|=\|M_1\|_{op} + \|M_2\|_{op}$, where $\|M_1\|_{op}$ and $\|M_2\|_{op}$ are the operator norms of the complex matrices $M_1$ and $M_2$ acting on $\mathbb{C}_i^{d\times1}$, respectively. Then we have the following theorem.


\begin{theorem} \label{bicomplex}
Let $d\in\mathbb{N}$, $\Lambda\subset[0,\infty)$ be an additive semigroup with $0\in\Lambda$, and let $\omega$ be an admissible weight on $\Lambda$. A function $f\in\mathcal{D}(\Lambda,\omega,M_d(\mathbb{BC}))$ has an inverse in $\mathcal{D}(\Lambda,\omega,M_d(\mathbb{BC}))$ if and only if $\overline{\ima(f)}=\overline{\{f(s):s\in\mathcal{H}\}}\subset GL_d(\mathbb{BC})$.
\end{theorem}

Lastly, we take $\mathcal{A}$ to be $M_d(\mathbb{H})$ ($d\in\mathbb{N}$), the real algebra of $d\times d$ matrices having entries from the real algebra of quaternions $\mathbb{H}$. We give brief details about the algebra of quaternions. An element $p\in\mathbb{H}$ is of the form $p=x_0+x_1 e_1+x_2 e_2+x_3 e_3$, where $x_0,x_1,x_2,x_3$ are real numbers, $e_1^2=e_2^2=e_3^2=e_1e_2e_3=-1$ and its norm is given by $|p|=\sqrt{x_0^2+x_1^2+x_2^2+x_3^2}$. Let $\mathbb{S}$ be the sphere of unitary purely imaginary quaternions, that is, it contains all $p\in\mathbb{H}$ such that $p^2=-1$. Two elements $i,j$ of $\mathbb{S}$ are \emph{\textbf{orthogonal}} if $ij+ji=0$ and they form a new basis $\{i,j,ij\}$ of $\mathbb{H}$. A matrix $M\in M_d(\mathbb{H})$ can then be uniquely written as $M=Z+Wj$, where $Z$ and $W$ are from $M_d(\mathbb{C}_i)$ and the norm of $M$ will be operator norm as an operator acting on $\mathbb{H}^{d\times 1}$. Then we have the following result. 

\begin{theorem} \label{quaternion}
Let $d\in\mathbb{N}$, $\Lambda\subset[0,\infty)$ be an additive semigroup with $0\in\Lambda$, and let $\omega$ be an admissible weight on $\Lambda$. A function $f\in\mathcal{D}(\Lambda,\omega,M_d(\mathbb{H}))$ has an inverse in $\mathcal{D}(\Lambda,\omega,M_d(\mathbb{H}))$ if and only if $\overline{\ima(f)}=\overline{\{f(s):s\in\mathcal{H}\}}\subset GL_d(\mathbb{H})$.
\end{theorem}

\begin{remark}
It is worth noting that in case of the bicomplex numbers we consider the \emph{commutative} orthogonal elements $i,j$ for generating a basis over $\mathbb{C}_i$ while in case of quaternion $i,j$ are chosen to be \emph{orthogonal} and this is the main difference giving different algebraic structures on these algebras. In fact in case of quaternion, $Z+Wj\neq Z+jW$ for $Z,W\in\mathbb{C}_i$.
\end{remark}

From Theorem \ref{matrix}, Theorem \ref{bicomplex} and Theorem \ref{quaternion}, we have the following corollary.

\begin{corollary}
Let $\Lambda\subset[0,\infty)$ be an additive semigroup with $0\in\Lambda$, and let $\omega$ be an admissible weight on $\Lambda$. For $d\in\mathbb{N}$, let $\mathcal{A}$ be one of the algebra $M_d(\mathbb{C})$, $M_d(\mathbb{BC})$ or $M_d(\mathbb{H})$. Then $\mathcal{D}(\Lambda,\omega,\mathcal{A})$ is inverse-closed in $\mathcal{D}(\Lambda,\mathcal{A})$.
\end{corollary}

It shall be noted that the results stated here trace back to the classical Wiener's theorem \cite{wi} for Fourier series and its weighted version obtained in \cite{bh,do} and references there in. The vector valued weighted versions are obtained in \cite{bo,kb}. Wiener's theorem for bicomplex numbers and quaternion valued functions are studied in \cite{bc} and \cite{al} respectively and they have inspired the results stated here for Dirichlet series. There are several applications of the Wiener's theorem in different fields of analysis like sampling theory, signal theory, Gabor analysis, frame theory, pseudo-differential analysis, Integral operator theory and many more; we refer the reader to \cite{Badi,Ba,kb3,KG,KB} and references therein. Of course the list is not exhaustive.

The proofs are provided in Section \ref{sec:proofs} and it is organized as follows: the proof for commutative algebra with $\Lambda=\log\mathbb{N}$ is provided in Section \ref{ss:commutative} and for complex, bicomplex and quaternion matrices in Section \ref{ss:complexmatrices}, \ref{ss:bicomplexmatrices} and \ref{ss:quaternionmatrices}, respectively.

\section{Proofs} \label{sec:proofs}
\subsection{Commutative algebra} \label{ss:commutative}
Let $\mathcal{A}$ be a commutative complex Banach algebra. A nonzero linear map $\varphi:\mathcal A \to \mathbb C$ satisfying $\varphi(ab)=\varphi(a)\varphi(b)\;(a,b \in \mathcal A)$ is a \emph{\textbf{complex homomorphism}} on $\mathcal A$. Let $\Delta(\mathcal A)$ be the collection of all complex homomorphisms on $\mathcal A$. For $a \in \mathcal A$, let $\widehat a:\Delta(\mathcal A)\to \mathbb C$ be $\widehat a(\varphi)=\varphi(a)\;(\varphi\in \Delta(\mathcal A))$. The smallest topology on $\Delta(\mathcal A)$ making each $\widehat a$, $a\in \mathcal A$, continuous is the \emph{\textbf{Gel'fand topology}} on $\Delta(\mathcal A)$ and $\Delta(\mathcal A)$ with the Gel'fand topology is the \emph{\textbf{Gel'fand space}} of $\mathcal{A}$. For more details on it refer \cite{bd}. We state a lemma which will be used in the proofs of Theorem \ref{commutative} and Theorem \ref{weighted commutative}.

\begin{lemma}\cite[Theorem 4 (v), pp. 82]{bd} \label{bdlem}
Let $\mathcal{A}$ be a commutative unital complex Banach algebra. An element $a\in\mathcal{A}$ is invertible if and only if $\phi(a)\neq0$ for all $\phi\in\Delta(\mathcal{A})$.
\end{lemma}

Considering $\mathbb{N}$ as a semigroup with multiplication, a nonzero bounded map $\chi:\mathbb N \to \mathbb C$ is a \emph{\textbf{semicharacter}} if $\chi(mn)=\chi(m)\chi(n)$ for all $m,n \in \mathbb N$. Let $\widehat{\mathbb{N}}$ be the collection of all semicharacters on $\mathbb{N}$. Then $\widehat{\mathbb{N}}$ is isomorphic the Cartesian product of countably infinite copies of $\mathbb{D}$, where $\mathbb D=\{z \in \mathbb C:|z|\leq 1\}$. For more details on it refer to \cite{he}. First we state a lemma which will be used for proving Theorem \ref{commutative}.

\begin{lemma}\label{l11}
Let $\mathcal{A}$ be a commutative unital complex Banach algebra. An element $a=(a_n)\in\ell^1(\mathbb{N},\mathcal{A})$ is invertible in $\ell^1(\mathbb{N},\mathcal{A})$ if and only if $\widehat{a}(\chi)=\sum_{n\in\mathbb{N}} a_n\chi(n)$ is invertible in $\mathcal{A}$ for all $\chi\in\widehat{\mathbb{N}}$.
\end{lemma}
\begin{proof}
It follows from \cite[Lemma 1]{he} and \cite[Theorem 3]{bo}.
\end{proof}

Following is the proof of Theorem \ref{commutative} which is inspired by the proof of Hewitt and Williamson of \cite[Theorem 1]{he}.

\begin{proof}[Proof of Theorem \ref{commutative}]
The equivalence of (i) and (iii) follows immediately. For (iii) implies (ii), let $x\in\overline{\ima(f_a)}$. Then there is a sequence $\{s_n\}_{n\in\mathbb{N}}$ in $\mathcal{H}$ such that $\lim_{n\to\infty}f_a(s_n)=x$. For $m,n\in\mathbb{N}$, \begin{align*} \|f_a(s_m)^{-1} - f_a(s_n)^{-1}\| &= \|f_a(s_n)^{-1} (f_a(s_n) - f_a(s_m)) f_a(s_m)^{-1} \| \\ &\leq \|f_a(s_n)^{-1}\| \|f_a(s_n) - f_a(s_m)\| \|f_a(s_m)^{-1} \| \\ &\leq \|b\|^2 \|f_a(s_n) - f_a(s_m)\| 
\end{align*} 
as $\displaystyle \|f_a(s)^{-1}\| = \left\|\sum_{n\in \mathbb N} b_n n^{-s} \right\| \leq \sum_{n\in \mathbb N}\|b_n\| = \|b\|$ for all $s\in\mathcal{H}$. This implies that $\{f_a(s_n)^{-1}\}_{n\in\mathbb{N}}$ is a Cauchy sequence in $\mathcal{A}$ and thus there is some $y\in\mathcal{A}$ such that $\lim_{n\to\infty}f_a(s_n)^{-1}=y$. This along with $\lim_{n\to\infty}f_a(s_n)=x$ gives $xy=1_\mathcal{A}$, that is, $x\in G(\mathcal{A})$. Thus, $\overline{\ima(f_a)}\subset G(\mathcal{A})$.

For (ii) implies (i), let $a\in\ell^1(\mathbb{N},\mathcal{A})$ be a singular element. Then, by Lemma \ref{l11}, there is a $\chi\in\widehat{\mathbb{N}}$ such that $\widehat{a}(\chi)=\sum_{n\in\mathbb{N}}a_n \chi(n)$ is not invertible in $\mathcal{A}$, that is, $\widehat{a}(\chi)\notin G(\mathcal{A})$. By Lemma \ref{bdlem}, there is some $\phi\in\Delta(\mathcal{A})$, the Gel'fand space of $\mathcal{A}$, such that $\phi(\widehat{a}(\chi))=\sum_{n\in\mathbb{N}}\phi(a_n)\chi(n)=0$. Let $\epsilon>0$. Then by techniques used in \cite[Theorem 1]{he}, there is a $j\geq0$ and a $t\in\mathbb{R}$ such that $|\phi(f_a(j+it))|=|\sum_{n\in\mathbb{N}} \phi(a_n) n^{-j-it}|<3\epsilon$. Since $\epsilon$ was arbitrary, it follows that $|\phi(f_a(s))|$ is not bounded away from zero on $\mathcal H$. This implies that $0\in\overline{\{\phi(f_a(s)):s\in\mathcal{H}\}}$. Since $\phi$ is continuous, it is equivalent to the fact that there is some $x\in\overline{\ima(f_a)}$ such that $\phi(x)=0$. This proves the theorem.
\end{proof}

Next is a proof of Theorem \ref{weighted commutative}, the weighted version Theorem \ref{commutative}.

\begin{proof}[Proof of Theorem \ref{weighted commutative}]
Since $\omega$ is an admissible weight, the Gel'fand space of $\mathcal{D}(\mathbb{N},\omega,\mathcal{A})$ is same as the Gel'fand space of $\mathcal{D}(\mathbb{N},\mathcal{A})$. By Lemma \ref{bdlem}, an element $f$ of $\mathcal{D}(\mathbb{N},\omega,\mathcal{A})$ is invertible in $\mathcal{D}(\mathbb{N},\omega,\mathcal{A})$ if and only if it is invertible in $\mathcal{D}(\mathbb{N},\mathcal{A})$ if and only if $\overline{\ima(f)} \subset G(\mathcal{A})$ by Theorem \ref{commutative}.
\end{proof}

\subsection{Complex matrices} \label{ss:complexmatrices}
First we state a lemma which establishes the relation between the scalar and matrix valued cases.

\begin{lemma}\label{equiv}
For an element $f\in\mathcal{D}(\Lambda,M_d(\mathbb{C}))$, the following are equivalent.
\begin{enumerate}
	\item $\overline{\ima(f)}=\overline{\{f(s):s\in\mathcal{H}\}}\subset GL_d(\mathbb{C})$.
	\item There is some $\delta>0$ such that for all $s\in\mathcal{H}$, $$|\det(f(s))|=\left|\det \begin{pmatrix} \sum_{\lambda\in\Lambda}f_{\lambda_{11}}e^{-\lambda s} & \sum_{\lambda\in\Lambda}f_{\lambda_{12}}e^{-\lambda s} & \dots & \sum_{\lambda\in\Lambda}f_{\lambda_{1d}}e^{-\lambda s} \\ \sum_{\lambda\in\Lambda}f_{\lambda_{21}}e^{-\lambda s} & \sum_{\lambda\in\Lambda}f_{\lambda_{22}}e^{-\lambda s} & \dots & \sum_{\lambda\in\Lambda}f_{\lambda_{2d}}e^{-\lambda s} \\ \vdots & \vdots &  \vdots & \vdots \\ \sum_{\lambda\in\Lambda}f_{\lambda_{d1}}e^{-\lambda s} & \sum_{\lambda\in\Lambda}f_{\lambda_{d2}}e^{-\lambda s} & \dots & \sum_{\lambda\in\Lambda}f_{\lambda_{dd}}e^{-\lambda s} \end{pmatrix}\right|\geq\delta,$$ where $f_{\lambda_{ij}}$ denotes the $ij^{th}$ entry of the matrix $f_\lambda$ for $i,j\in\{1,2,\dots,d\}$.
\end{enumerate}
\end{lemma}

\begin{proof}
(ii) is not true if and only if there is a sequence $\{s_n\}_{n\in\mathbb{N}}$ in $\mathcal H$ such that $\lim_{n\to\infty}\det(f(s_n))=0$ if and only if $\lim_{n\to\infty} f(s_n)\notin GL_d(\mathbb{C})$ if and only if (i) is not true.
\end{proof}

So, when we have $d=1$, then the above lemma shows that the condition mentioned in Theorem \ref{matrix} is same as the one taken in Theorem \ref{th1} and \ref{th2}. Also, observe that for a matrix $A=(a_{kl})\in M_d(\mathbb C)$, let $\|A\|_1=\sum_{k,l=1}^n |a_{kl}|$. Then $\|\cdot\|_1$ is an algebra norm on $M_d(\mathbb C)$ and since $M_d(\mathbb C)$ is finite dimensional, the operator norm $\|\cdot\|_{op}$ and the one norm $\|\cdot\|_1$ are equivalent. Now we give a proof of Theorem \ref{matrix}.

\begin{proof}[Proof of Theorem \ref{matrix}]
Let $f\in\mathcal{D}=\mathcal{D}(\Lambda,\omega,M_d(\mathbb{C}))$ be invertible in $\mathcal{D}$. Then there is some $g\in\mathcal{D}$ such that $fg=I_d$. This means that $(\sum_{\lambda\in\Lambda} f_\lambda e^{-\lambda s})(\sum_{\lambda\in\Lambda} g_\lambda e^{-\lambda s}) = I_d$ for all $s\in\mathcal{H}$. Then $\overline{\ima(f)}\subset GL_d(\mathbb{C})$ follows from same arguments used in the proof of (iii) implies (ii) of Theorem \ref{commutative} by taking $\mathcal{A}=M_d(\mathbb{C})$ as the commutativity of $\mathcal{A}$ is not used anywhere there.

Conversely, assume that $f\in\mathcal{D}$ is such that $\overline{\ima(f)}=\overline{\{f(s):s\in\mathcal{H}\}}\subset GL_d(\mathbb{C})$. Then for given $s\in\mathcal{H}$, there is a matrix $G_s$ such that $f(s)G_s=I_d$. This matrices $G_s$ are given by 
$$G_s=\frac{1}{\det(f(s))}\begin{pmatrix} A_{11} & A_{12} & \dots & A_{1d} \\ \vdots & \vdots & & \vdots \\ A_{d1} & A_{d2} & \dots & A_{dd} \end{pmatrix},$$ 
where $A_{kl}$, for $k,l\in\{1,2,\dots,d\}$, are the cofactors of the matrix 
$$f(s)=\begin{pmatrix} \sum_{\lambda\in\Lambda}f_{\lambda_{11}}e^{-\lambda s} & \sum_{\lambda\in\Lambda}f_{\lambda_{12}}e^{-\lambda s} & \dots & \sum_{\lambda\in\Lambda}f_{\lambda_{1d}}e^{-\lambda s} \\ \sum_{\lambda\in\Lambda}f_{\lambda_{21}}e^{-\lambda s} & \sum_{\lambda\in\Lambda}f_{\lambda_{22}}e^{-\lambda s} & \dots & \sum_{\lambda\in\Lambda}f_{\lambda_{2d}}e^{-\lambda s} \\ \vdots &\vdots & \vdots & \vdots \\ \sum_{\lambda\in\Lambda}f_{\lambda_{d1}}e^{-\lambda s} & \sum_{\lambda\in\Lambda}f_{\lambda_{d2}}e^{-\lambda s} & \dots & \sum_{\lambda\in\Lambda}f_{\lambda_{dd}}e^{-\lambda s} \end{pmatrix}.$$
Since each entry of the matrix $f(s)$ is an element of $\mathcal{D}(\Lambda,\omega,\mathbb{C})$ and addition and multiplication of elements of $\mathcal{D}(\Lambda,\omega,\mathbb{C})$ are again in $\mathcal{D}(\Lambda,\omega,\mathbb{C})$, each $A_{kl}$ is in $\mathcal{D}(\Lambda,\omega,\mathbb{C})$ for $k,l\in\{1,2,\dots,d\}$. Also, $\det(f(s))$ is an element $\mathcal{D}(\Lambda,\omega,\mathbb{C})$ as $\det(f(s))$ is obtained by adding and multiplying the elements of the Banach algebra $\mathcal{D}(\Lambda,\omega,\mathbb{C})$. Since $\overline{\ima(f)}\subset GL_d(\mathbb{C})$, by Lemma \ref{equiv}, there is some $\delta>0$ such that $|\det(f(s))|\geq\delta$ for all $s\in\mathcal H$. That is $\det(f(s))$ is an element of $\mathcal{D}(\Lambda,\omega,\mathbb{C})$ which is bounded away from 0 in the absolute value and so, by Theorem \ref{th2}, it has an inverse in the form of a Dirichlet series, say $\sum_{\lambda\in\Lambda}h_\lambda e^{-\lambda s}$ in $\mathcal{D}(\Lambda,\omega,\mathbb{C})$. Thus, we have 
\begin{align*}G_s=\sum_{\lambda\in\Lambda}h_\lambda e^{-\lambda s} \begin{pmatrix} A_{11} \ A_{12} \ \dots \ A_{1d} \\ \vdots \\ A_{d1} \ A_{d2} \ \dots \ A_{dd} \end{pmatrix}. \end{align*} 
Again using the same fact as above, it follows that $G_s$ is a $d\times d$ matrix having each entry as an element of $\mathcal{D}(\Lambda,\omega,\mathbb{C})$. Define $g(s)=G_s$ $(s\in\mathcal{H})$. Then the uniqueness of Dirichlet series gives that $g$ is of the form $$g(s)=\sum_{\lambda\in\Lambda} g_\lambda e^{-\lambda s} \quad (s\in\mathcal{H}),$$ where $g_\lambda\in M_d(\mathbb{C})$ for all $\lambda\in\Lambda$. In order to show that $\|g\|_\omega=\sum_{\lambda\in\Lambda} \|g_\lambda\| \omega(\lambda)<\infty$, note that 
$$g(s)=\begin{pmatrix} \sum_{\lambda\in\Lambda} g_{\lambda_{11}}e^{-\lambda s} & \sum_{\lambda\in\Lambda} g_{\lambda_{12}}e^{-\lambda s} & \dots & \sum_{\lambda\in\Lambda} g_{\lambda_{1d}}e^{-\lambda s} 
\\ \sum_{\lambda\in\Lambda} g_{\lambda_{21}}e^{-\lambda s} & \sum_{\lambda\in\Lambda} g_{\lambda_{22}}e^{-\lambda s} & \dots & \sum_{\lambda\in\Lambda} g_{\lambda_{2d}}e^{-\lambda s} \\ \vdots  & & \\ \sum_{\lambda\in\Lambda} g_{\lambda_{d1}}e^{-\lambda s} & \sum_{\lambda\in\Lambda}g_{\lambda_{d2}}e^{-\lambda s} & \dots & \sum_{\lambda\in\Lambda} g_{\lambda_{dd}}e^{-\lambda s} \end{pmatrix},$$ 
where each $g_{\lambda_{kl}}$ for $k,l\in{1,2,\dots,d}$ and $\lambda\in\Lambda$ is obtained as multiplication of elements of $\mathcal{D}(\Lambda,\omega,\mathbb{C})$ and so $\sum_{\lambda\in\Lambda} |g_{\lambda_{kl}}| \omega(\lambda)  <\infty$ for all $k,l\in{1,2,\dots,d}$. Thus, $$\sum_{\lambda\in\Lambda} \sum_{k,l=1}^d |g_{\lambda_{kl}}| \omega(\lambda)=\sum_{\lambda\in\Lambda} \|g_\lambda\|_1 \omega(\lambda)<\infty.$$ Since the one norm and operator norm are equivalent on $M_d(\mathbb C)$, it follows that $\displaystyle \|g\|_\omega=\sum_{\lambda\in\Lambda} \|g_\lambda\| \omega(\lambda)<\infty.$ This completes the proof.
\end{proof}

\subsection{Bicomplex matrices} \label{ss:bicomplexmatrices}

\begin{proof}[Proof of Theorem \ref{bicomplex}]
Let $f\in\mathcal{D}=\mathcal{D}(\Lambda,\omega,M_d(\mathbb{BC}))$ be such that $\overline{\ima(f)}\subset GL_d(\mathbb{BC})$. Then $f$ can be written as $$f(s)=\sum_{\lambda\in\Lambda} (f_{\lambda1}\mathbf{e}_1 e^{-\lambda s} + f_{\lambda2}\mathbf{e}_2 e^{-\lambda s}).$$ Since $\overline{\ima(f)}\subset GL_d(\mathbb{BC})$ if and only if $\ima(\sum_{\lambda\in\Lambda} f_{\lambda1} e^{-\lambda s})$ and $\ima(\sum_{\lambda\in\Lambda} f_{\lambda2} e^{-\lambda s})$ are contained in $GL_d(\mathbb{C})$, by Theorem \ref{matrix}, there are $\sum_{\lambda\in\Lambda} g_{\lambda k} e^{-\lambda s}$ in $\mathcal{D}(\Lambda,\omega,M_d(\mathbb{C}))$ which are inverses of $\sum_{\lambda\in\Lambda} f_{\lambda k} e^{-\lambda s}$ for $k=1,2$. Define $$g(s)=\sum_{\lambda\in\Lambda} (g_{\lambda1} \mathbf{e}_1 e^{-\lambda s} + g_{\lambda2} \mathbf{e}_2 e^{-\lambda s}).$$ Then $g$ is our required function. The converse follows from the arguments similar to the one given in Theorem \ref{matrix}.
\end{proof}

\subsection{Quaternion matrices} \label{ss:quaternionmatrices}
Before moving on to the proof of Theorem \ref{quaternion}, we have the following lemma whose proof is straightforward and hence not included. For a matrix $M=(m_{kl})\in M_d(\mathbb{C})$, the matrix $\overline{M}=(\overline{m_{kl}})$ is again in $M_d(\mathbb{C})$. Here, $\overline{z}$ denotes the complex conjugate of $z\in\mathbb{C}$.

\begin{lemma}
If $d\in\mathbb{N}$, $\Lambda\subset[0,\infty)$ is an additive semigroup with $0\in\Lambda$, $\omega$ is an admissible weight on $\Lambda$, and if $f(s)=\sum_{\lambda\in\Lambda} f_\lambda e^{-\lambda s}$ $(s\in\mathcal{H})$ is in $\mathcal{D}(\Lambda,\omega,M_d(\mathbb{C}))$, then $\overline{f}(s)=\sum_{\lambda\in\Lambda} \overline{f_\lambda} e^{-\lambda s}$ is again in $\mathcal{D}(\Lambda,\omega,M_d(\mathbb{C}))$. Moreover, if $f$ has an inverse in $\mathcal{D}(\Lambda,\omega,M_d(\mathbb{C}))$, then $\overline{f}$ also has an inverse in $\mathcal{D}(\Lambda,\omega,M_d(\mathbb{C}))$.    
\end{lemma}

Now we are ready for the proof of Theorem \ref{quaternion}.

\begin{proof}[Proof of Theorem \ref{quaternion}]
We give the proof for the case of $d=1$. Let $j,l\in\mathbb{S}$ be orthogonal quaternions. Then any quaternion $h\in\mathbb{H}$ can be written as $h=h^j+h^l l$ where $h^j,h^l$ are from the complex plane $\mathbb{C}_j=\{x+jy:x,y\in\mathbb{R}\}$. Let $f\in\mathcal{D}(\Lambda,\omega,\mathbb{H})$. Then we may write $f$ as \begin{align*}f(s)&=\left( \sum_{\lambda\in\Lambda} f^j_\lambda(s) + \sum_{\lambda\in\Lambda} f^l_\lambda(s) l \right) e^{-\lambda s} \\&=\sum_{\lambda\in\Lambda} f^j_\lambda(s) e^{-\lambda s} + \sum_{\lambda\in\Lambda} f^l_\lambda(s) l e^{-\lambda s}  \\&= f^j(s) + f^l(s) l  \quad (s\in\mathcal{H}).\end{align*} Suppose that $g(s)=\left(\sum_{\lambda\in\Lambda} g^j(s) + \sum_{\lambda\in\Lambda} g^l(s) l \right) e^{-\lambda s}$ $(s\in\mathcal{H})$ is an inverse of $f$. Then for each $s\in\mathcal{H}$ \begin{align*} &(f^j(s) + f^l(s)l) (g^j(s) + g^l(s)l) = 1 \\ \implies &f^j(s)g^j(s) + f^j(s) g^l(s) + f^l(s)lg^j(s) + f^l(s)l g^l(s)l=1 \\ \implies &f^j(s)g^j(s) - f^l(s) \overline{g^l}(s) + [f^j(s)g^l(s) + f^l(s) \overline{g^j}(s)]l=1 \\ \implies &f^j(s)g^j(s) - f^l(s) \overline{g^l}(s)=1 \ \text{and} \ f^j(s)g^l(s) + f^l(s) \overline{g^j}(s)=0. \end{align*} Simplifying these two equation, it follows that $$g^j(s)=\frac{\overline{f^j}(s)}{|f^j(s)|^2 + |f^l(s)|^2} \ \text{and} \ g^l(s)=\frac{-f^l(s)}{|f^j(s)|^2 + |f^l(s)|^2} .$$ Since $g\in\mathcal{D}(\Lambda,\omega,\mathbb{H})$, $g^j,g^l\in\mathcal{D}(\Lambda,\omega,\mathbb{C})$ and thus there is some $\delta>0$ such that $|f^j(s)|^2 + |f^l(s)|^2\geq\delta$ for all $s\in\mathcal{H}$, otherwise the functions $g^j,g^l$ cannot be in $\mathcal{D}(\Lambda,\omega,\mathbb{C})$. It is equivalent to the fact that $\overline{\ima(f)}\subset\mathbb{H}\setminus\{0\}$. The converse follows by defining $g$ as above and using Theorem \ref{th2} for $g^j$ and $g^l$.

The case of $d>1$ can be dealt with proper modifications as done in Theorem \ref{matrix} and Theorem \ref{bicomplex}. 
\end{proof}

\subsection*{Further remarks}
\begin{remark}
It shall be noted that for $\mathbb{X}\in\{\mathbb{C},\mathbb{BC},\mathbb{H}\}$, the case where $\mathcal{A}$ is $\mathbb{X}^d$ for arbitrary $\Lambda$ and $d\in\mathbb{N}$ follows by replacing $M_d(\mathbb{X})$ by the collection of all diagonal matrices having entries from $\mathbb{X}$ and using the fact that this collection is isomorphic to $\mathbb{X}^d$.
\end{remark}

\begin{remark}
We do not know whether the results stated above follows for arbitrary algebra valued Dirichlet series. That is we have the following open question.
\begin{itemize}
    \item [] Let $\Lambda\subset[0,\infty)$ be an additive semigroup with $0\in\Lambda$, $\omega$ be an admissible weight on $\Lambda$, $\mathcal{A}$ be a complex or real Banach algebra (not necessarily commutative), and let $f\in\mathcal{D}(\Lambda,\omega,\mathcal{A})$ is such that $\overline{\ima(f)}\subset G(\mathcal{A})$. Is $f$ invertible in $\mathcal{D}(\Lambda,\omega,\mathcal{A})$?
\end{itemize}
\end{remark}

\begin{remark}
The difficulty in dealing with arbitrary $\Lambda$ in place of $\log\mathbb{N}$ for commutative complex Banach algebra case lies in the fact that the proof provided here relies on the Gel'fand theory and even the case of $\Lambda=\mathbb{Q}_+$, set of non-negative rationals, can not be dealt easily because the Gel'fand space of $\ell^1(\mathbb{Q}_+)$ is not known explicitly. The study of identifying the Gel'fand spaces of different algebras in explicit form is of importance on its own and not easy.
\end{remark}

\section*{Statements and Declarations}
\textbf{Funding.} This work is not funded by any specific project. The second author gratefully acknowledges the financial support from Post Doctoral Fellowship under the ISIRD project 9--551/2023/IITRPR/10229 of Dr. Manmohan Vashisth from IIT Ropar.

\textbf{Competing Interests.} The authors have no relevant financial or non-financial interests to disclose.

\textbf{Author Contributions.} All authors contributed equally.

\textbf{Availability of data and material.} Not applicable.

\bibliographystyle{amsplain}

\end{document}